\documentclass[12pt]{article}
\usepackage{amssymb,theorem,amsmath,enumerate}
\usepackage{graphicx}
\usepackage{tikz}
\usepackage{slashbox}
\usetikzlibrary{arrows,decorations.pathmorphing,backgrounds,positioning,fit,petri}
\newcommand{\bc}{\begin{center}}
\newcommand{\ec}{\end{center}}
\newcommand{\be}{\begin{enumerate}}
\newcommand{\ee}{\end{enumerate}}
\newcommand{\beq}{\begin{equation}}
\newcommand{\eeq}{\end{equation}}

\newcommand{\ts}{\textstyle}
\newcommand{\lb}{\linebreak}
\newcommand{\ol}{\overline}
\newcommand{\bb}{\rule{1.5mm}{1.5mm}}
\newcommand{\ba}{\begin{array}}
\newcommand{\ea}{\end{array}}
\newcommand{\bea}{\begin{eqnarray*}}
\newcommand{\eea}{\end{eqnarray*}}
\newcommand{\mb}{\mbox}
\newcommand{\Z}{\Bbb Z}
\newcommand{\N}{\Bbb N}
\newcommand{\Q}{\Bbb Q}
\newcommand{\ra}{\rightarrow}

\begin{document}
\renewcommand{\thefootnote}{\fnsymbol{footnote}}
\setcounter{MaxMatrixCols}{20}
\include{BEcommandsExpInlineVar}

\bc {\bf \Large Functions realising as abelian group automorphisms} \ec

\bc {\bf B-E de Klerk, JH Meyer, J Szigeti, L van Wyk} \ec

{\bf Abstract:} Let $A$ be a set and  $f:A\ra A$ a bijective function. Necessary and sufficient conditions on $f$ are determined which makes it possible to endow $A$ with a binary operation $*$ such that $(A,*)$ is a cyclic group and $f\in \mb{Aut}(A)$. This result is extended to all abelian groups in case $|A|=p^2, \ p$ a prime. Finally, in case $A$ is countably infinite, those $f$ for which it is possible to turn $A$ into a group $(A,*)$ isomorphic to $\Z^n$ for some $n\ge 1$, and with $f\in \mb{Aut} (A)$, are completely characterised.  \\

{\sl Keywords:} Automorphism, abelian group

{\sl 2010 Mathematics Subject Classification:} 20K30, 20K01, 20E34\\

{\bf 1. Introduction}

The question on which functions from a set to itself (selfmaps) appear as functions with a certain structural property, has been addressed by various authors. In particular, in \cite{FolSzi} and \cite{Szigeti} those selfmaps which appear as lattice endomorphisms or lattice anti-endomorphisms have been characterised. In \cite{BouFonKeYeh} a similar study was done for infra-endomorphisms of the groups $\Z_n$ and $D_n$. In this paper we characterise those selfmaps that appear as automorphisms of certain abelian groups, namely the cyclic groups, the group $\Z_p\times \Z_p$, $p$ prime, and the group $\Z^n$ for some $n\ge1$.  

For a given set $A$, let us agree to say that a bijection $f:A\ra A$ has the {\it auto-property} if it is possible to find a binary operation $*$ on $A$ such that $(A,*)$ is an abelian group and $f\in \mb{Aut}(A)$. 

If $A$ is finite, such an $f$ necessarily gives rise to cycles, i.e., (disjoint) finite sequences $a_1,a_2,\ldots,a_m$ from $A$ such that $f(a_i)=a_{i+1}$ for $1\le i \le m-1$ and $f(a_m)=a_1$. Every element of $A$ belongs to some cycle. The number of elements in a cycle is its {\it length}. So a fixed point of $f$ is a cycle of length $1$. The {\it cycle structure} of $f$ is a description of how many cycles of each length $f$ has. A convenient notation for this structure will be developed and used in Section 2.  

On the other hand, if $A$ is infinite, then, apart from possible cycles, there is also the possibility of $f$ having {\it chains}, i.e., infinite sequences $\ldots, a_i, a_{i+1}, a_{i+2},\ldots$ from $A$ such that $f(a_i) = a_{i+1}$ for all $i$. The number of cycles of various lengths, as well as the number of chains, will be referred to simply as the {\it structure} of $f$.  This infinite case will be discussed in Section 3.  \\

{\bf 2. Cyclic groups and groups of order $p^2$, $p$ prime}

For this section, we always assume that $A$ is a finite set and that $f:A\ra A$ is a bijection.
 
If $f$ has $c_i$ cycles of length $t_i$ ($1\le i \le k$), then we say $f$ has the {\it cycle structure} $\left[\ba{cccc}c_1&c_2&\cdots&c_k\\t_1&t_2&\cdots&t_k\ea\right]$, where, for consistency, we always take $t_1>t_2>\cdots>t_k$. We stress that it is possible that some $c_i$ could be $0$, and these columns can just as well be omitted from the array.  
Also note that $\sum_{i=1}^k c_it_i=|A|$, and that the identity map has cycle structure $\left[\ba{c}|A|\\1\ea\right]$. 

It is evident that if we want to investigate the conditions $f$ has to satisfy to have the auto-property, then it suffices to find the cycle structures of all possible automorphisms on all abelian groups of order $|A|$. 
These cycle structures completely determine all those $f$ having the auto-property. 

We begin by doing this for cyclic groups. 

We determine all possible cycle structures of automorphisms $f:\Z_n \ra \Z_n$, for the additive (cyclic) group $\Z_n = \{0,1,\ldots,n-1\}$.

For the ring $\Z_n = \{0,1,\ldots,n-1\}$, consider the group of units $U_n = \{k_1,k_2,\ldots,k_{\phi(n)}\} = \{k \in \Z_n : (k,n)=1\}$ (where we take $k_1 = 1$).  Let $T_n = (\Z_n\setminus U_n)\setminus \{0\}$ and  
for $z\in T_n$, put $z' = \frac{n}{(z,n)}$. 

Let $f:\Z_n \ra \Z_n$ be an automorphism. Then $f(1)\in U_n$, otherwise, if $f(1) = z\in T_n$, then $f(z') = 0 = f(0)$, a contradiction. If $f(1) = 1 = k_1$, then $f$ is the identity map.

 Let $2\le i \le \phi(n)$, and assume that $f(1) = k_i$.  Then $1, k_i, k_i^2,\ldots,k_i^{\ell_i-1}$ is a cycle of length $\ell_i = \mb{ord}_n(k_i)$ 
(the least $x\in\N$ such that $k_i^x\equiv 1(\mb{mod } n)$), and consisting exactly of the elements of the subgroup $\langle k_i \rangle$ of $U_n$. If $\langle k_i \rangle \ne U_n$, choose any $k_j \in U_n\setminus \langle k_i \rangle$. Then  $k_j, k_jk_i, k_jk_i^2,\ldots,k_jk_i^{\ell_i-1}$ is another cycle of length $\ell_i$, consisting exactly of the coset $k_j\langle k_i \rangle$ of $\langle k_i \rangle$ in $U_n$. Continuing in this manner, we obtain $[U_n : \langle k_i \rangle]$ cycles of this type, exhausting all the elements of $U_n$.
 
  Now consider any $z \in T_n$. Then the cycle $z, zk_i, zk_i^2,\ldots,zk_i^{\lambda-1}$ is obtained, where the length of the cycle is the least $\lambda\in\N$ such that $n\,|\,z(k_i^{\lambda}-1)$. This means that $\lambda=\mb{ord}_{z'}(k_i)$. Note that $\lambda|\ell_i$. Also note that each member of this cycle is in $T_n$. Other elements of $T_n$, not in this cycle, might give rise to cycles of the same length $\lambda$. Hence, the total number of cycles of length $\lambda$ is given by $\frac1{\lambda}\left|\{z\in T_n : \mb{ord}_{z'}(k_i)=\lambda\}\right|=:L_{i,\lambda}$. Finally, cycles of length $1$ obtained in this way exclude the fixed point $0$, so that there are $\left|\{z\in T_n : \mb{ord}_{z'}(k_i)=1\}\right|+1$ cycles of length $1$.
  
Hence we have one direction of the following theorem:

 {\bf 2.1. Theorem.\ } Let $|A|=n$ and let $f:A\ra A$ be a bijection. Then there exists a binary operation $*$ on $A$ such that $(A,*)$ is a cyclic group and $f\in \mb{Aut}(A)$ if and only if either $f$ is the identity map, or there is an $i,\ 2\le i \le \phi(n)$, such that $f$ has the cycle structure 
 $$\left[\ba{ccccc}[U_n:\langle k_i \rangle]+L_{i,\ell_i}&L_{i,\lambda_1}&\cdots&L_{i,\lambda_{t}}&L_{i,1} + 1\\\ell_i&\lambda_1&\cdots&\lambda_t&1\ea\right],$$ where $\ell_i > \lambda_1 > \cdots > \lambda_t>1$ denotes the complete list of (positive) divisors of $\ell_i=\mb{ord}_n(k_i)$. 
 
 {\bf Proof: } It remains to show how to turn $A$ into an abelian group with $f\in \mb{Aut}(A)$, given that $f$ satisfies the stated conditions . This can be done {\it via} the so-called {\it structural graph} of $f$. Let ${\cal G} = (V,E)$ be a directed graph with $|V| = n$ and $\rho : A\rightarrow V$ a bijection.  Then ${\cal G}$ is called a {\it structural graph} of $f$ if $(u,v) \in E \Leftrightarrow(\exists a\in A : u = \rho(a) \wedge v = \rho(f(a)))$. $\rho$ is called a {\it graph projection} of $f$. 
 
 Now, if there exists a group automorphism $h : G \rightarrow G$ for some abelian group $G$ such that the structural graphs of $f$ and $h$ are isomorphic (as graphs), then one easily sees that $A$ can be endowed with and abelian group structure such that $f$ is a group automorphism.
 
 In particular, let $\rho_f$ and $\rho_h$ be graph projections of $f$ and $h$ respectively, and let $\psi$ be a graph isomorphism from the codomain of $\rho_f$ to the codomain of $\rho_h$. Define $\eta : A\rightarrow G$ by $\eta = \rho_h^{-1}\psi\rho_f$. Then it is routine to check that $(A,*)$ is an abelian group, where $\alpha * \beta = \eta^{-1}(\eta(\alpha)\cdot_G \eta(\beta))$ for all $\alpha, \beta \in A$. The identity is $1_A = \eta^{-1}(1_G)$. It is also routine to check that $f\in \mb{Aut}(A)$.
   \ \bb
 

{\bf 2.2. Corollary.\ } If $|A|=n$, then there are at most $\phi(n)$ cycle structures for a bijection $f:A\ra A$ that will turn $A$ into a cyclic group, with $f\in\mb{Aut}(A)$.

{\bf  Proof: \ } Apart from the identity map, the possible cycle structures of automorphisms are determined by $2\le i \le \phi(n)$. But note that distinct $i$'s could give rise to the same cycle structure of an automorphism. \ \bb

{\bf 2.3. Example.\ } 
\be
\item[(a)] If $|A|=p$, where $p$ is a prime, then $f:A\ra A$ has the auto-property if and only if it has the cycle structure  $\left[\ba{cc}d&1\\\frac{p-1}{d}&1\ea\right]$ for some divisor $d$ of $p-1$. (Note that in case $d=p-1$, we get that $\left[\ba{cc}d&1\\\frac{p-1}{d}&1\ea\right] = \left[\ba{c}  p\\1\ea\right]$, representing the identity map.)

\item[(b)] Let $|A|=12$. Then $U_{12}=\{1,5,7,11\}$, so that $(k_1,k_2,k_3,k_4)=(1,5,7,11)$. Then we have 

 $\ell_2 = \mb{ord}_{12}(5)=2$. $L_{2,1}=|\{3,6,9\}|=3, \ L_{2,2}=\frac12\cdot|\{2,4,8,10\}|=2$. This gives the cycle structure $\left[\ba{cc} [U_{12}:\langle 5\rangle] +L_{2,2} & L_{2,1}+1\\ 2&1\ea\right] = \left[\ba{cc} 4 & 4\\ 2&1\ea\right]$. 
 
 Similarly, for $\ell_3 = 2$ we obtain the cycle structure $\left[\ba{cc} 3 & 6\\ 2&1\ea\right]$ and for $\ell_4 = \ell_{\phi(12)}= 2$ we obtain the cycle structure $\left[\ba{cc} 5 & 2\\ 2&1\ea\right]$. 
 
 Hence, $A$ can be turned into a cyclic group with $f\in \mb{Aut}(A)$ if and only if $f$ is the identity map, or $f$ has one of the three cycle structures above.
 
 \item[(c)] Let $|A|=p^2$, with $p$ prime. Then $z' = p$ for all $z\in T_{p^2}= \{p,2p,\ldots,(p-1)p\}$. This implies that
 \[ L_{i,\lambda} = \ts\frac1{\lambda}\left|\{z\in T_{p^2} : \mb{ord}_{p}(k_i)=\lambda\}\right| = \left\{\ba{cl } \frac{p-1}{\lambda} & \mb{if } \mb{ord}_p(k_i)=\lambda\\ 0 & \mb{otherwise}\ea\right.\]
 for every divisor $\lambda$ of $\ell_i = \mb{ord}_{p^2}(k_i)$, where $2\le i \le p^2-p$. 
 
 For instance, if $p=3$, then $(k_1,k_2,\ldots,k_6) = (1,2,4,5,7,8)$. For $i=2$ we have $\ell_2 = \mb{ord}_9(2) = 6$, and since $\mb{ord}_3(2) = 2$, it follows that $L_{2,2} = \frac22 = 1$ and $L_{2,1} = L_{2,3} = L_{2,6} = 0$. Also, since $k_2 = 2$ is a generator of the group $U_9$, $[U_9 : \langle 2 \rangle] = [U_9:U_9] = 1$. So (for the case $i=2$) we obtain, by Theorem 2.1, the cycle structure $\left[\ba{cccc} 1 &0& 1&1\\ 6&3&2&1\ea\right] =\left[\ba{ccc} 1 & 1&1\\ 6&2&1\ea\right]$. 
 
 Similarly, for $i=3$ we get the cycle structure $\left[\ba{cc} 2 & 3\\ 3&1\ea\right]$, for $i=4$ we get $\left[\ba{ccc} 1&1&1\\ 6&2&1\ea\right]$, for $i=5$ we get $\left[\ba{cc}  2&3\\ 3&1\ea\right]$, and finally, for $i=6=\phi(9)$ we get $\left[\ba{cc} 4&1\\ 2&1\ea\right]$. 
 
 Consequently, if $|A|=9$, it can be turned into a cyclic group with $f\in \mb{Aut}(A)$ if and only if $f$ has one of the cycle structures $\left[\ba{ccc} 1&1&1\\ 6&2&1\ea\right], \ \left[\ba{cc} 2 & 3\\ 3&1\ea\right], \ \left[\ba{cc} 4&1\\ 2&1\ea\right]$ or $\left[\ba{c} 9\\ 1\ea\right]$ (the identity). Note that other cycle structures are indeed possible in the non-cyclic case (see Theorem 2.5). 
  
 \ee
 
 We now turn our attention to the case $|A|=p^2$, $p$ prime, and completely determine when $f$ has the auto-property in this case. Theorem 2.1 takes care of the case when $A$ is cyclic. We will therefore focus here only on the group $\Z_p^2= \Z_p \times \Z_p$, with Aut$(\Z_p \times \Z_p) \cong GL_2(\Z_p)$.  Our aim is to determine the cycle structures of all the elements of $GL_2(\Z_p)$, when acting on the elements of $\Z_p^2$. 
 
We recall that conjugate  permutations have the same cycle structures, and we formalize this in

{\bf 2.4. Lemma.\ } If $F$ is a finite field, and $A,B\in GL_2(F)$ are similar, then they determine the same cycle structure on the group $F^2$. \ \bb   


Henceforth, for $\alpha$ in the finite field $F$,  we use the notation $o^+(\alpha)$ for the (additive) order of $\alpha\in F$ and we use $o^\times(\alpha)$ for the (multiplicative) order of $\alpha\in F^*$.

In \cite{glg} it is given that there exist elements of order $d$ in $GL_2(\Z_p)$, for any $d$ that divides  $p^2-1$, as well as of order $pd$ for any $d\,|\, p-1$.  Furthermore, by virtue of Lemma 2.4, we only have to study the Jordan normal forms of the matrices in $GL_2(\Z_p)$. We do it by considering three cases:

\be 
\item[I.] Here, we only consider those matrices $A$ in $GL_2(\Z_p)$ having Jordan normal form $\left(\ba{cc}\alpha_1 & 0\\0 & \alpha_2\ea\right)$, where $\alpha_1,\alpha_2\in U_p$. The order of such an $A$ is $d$, where $d\,|\,p-1$.  

First, if $\alpha_1=\alpha_2=\alpha$ (say), with $o^\times(\alpha)=d$, then $A=\left(\ba{cc}\alpha & 0\\0 & \alpha\ea\right)$ has $\frac{p^2-1}{d}$ cycles of the form $$\left(\ba{c} x\\y \ea\right), \left(\ba{c}\alpha x\\\alpha y \ea\right),\ldots, \left(\ba{c} \alpha^{d-1}x\\\alpha^{d-1}y \ea\right),$$ each of length $d$ and where $x,y\in\Z_p$, not both $0$.

Second,  let $\alpha_1\ne \alpha_2$, with $o^\times(\alpha_1)=d_1$ and $o^\times(\alpha_2)=d_2$, where $d_1$ and $d_2$ are divisors of $p-1$. Here,  $A=\left(\ba{cc}\alpha_1 & 0\\0 & \alpha_2\ea\right)$ has

$\frac{p-1}{d_1}$ cycles of the form
$$\left(\ba{c} x\\0 \ea\right), \left(\ba{c}\alpha_1 x\\0 \ea\right),\ldots, \left(\ba{c} \alpha_1^{d_1-1}x\\0 \ea\right),$$
each of length $d_1$, where $x\in U_p$;  

$\frac{p-1}{d_2}$ cycles of the form
$$\left(\ba{c} 0\\y \ea\right), \left(\ba{c}0 \\\alpha_2 y \ea\right),\ldots, \left(\ba{c} 0\\\alpha_2^{d_2-1}y \ea\right),$$
each of length $d_2$, where $y\in U_p$;

$\frac{(p-1)^2}{\mb{\footnotesize lcm}(d_1,d_2)}$ cycles of the form
$$\left(\ba{c} x\\y \ea\right), \left(\ba{c}\alpha_1 x \\\alpha_2 y \ea\right),\ldots,\left(\ba{c} \alpha_1^{K-1}x\\\alpha_2^{K-1}y \ea\right),$$
each of length $K=\mb{lcm}(d_1,d_2)$, where $x,y\in U_p$.

\item[II.] Now we consider those $A\in GL_2(\Z_p)$ with Jordan normal form $\left(\ba{cc}\alpha & 1\\0 & \alpha\ea\right)$, where $\alpha\in U_p$. The order of $A$ is $pd$, where $d=o^\times(\alpha)$ is a divisor of $p-1$, and for any $d\,|\,p-1$, there exists such an $A$.   

Then  $A$ has $\frac{p-1}{d}$ cycles of the form
$$\left(\ba{c} x\\0 \ea\right), \left(\ba{c}\alpha x\\0 \ea\right),\ldots, \left(\ba{c} \alpha^{d-1}x\\0 \ea\right),$$
each of length $d$, where $x\in U_p$;

 $A$ has $\frac{p-1}{d}$ cycles of the form
$$\left(\ba{c} x\\y \ea\right), \left(\ba{c}\alpha x+y \\ \alpha y \ea\right),\ldots,\left(\ba{c}\alpha^k x+k\alpha^{k-1} y\\ \alpha^k y \ea\right),\ldots, \left(\ba{c}\alpha^{pd-1}x+ (pd-1)\alpha^{pd-2}y\\\alpha^{pd-1}y \ea\right),$$
each of length $pd$, where $x,y\in \Z_p$, with $y\ne0$. (Note that $o^+(d\alpha^{pd-1}y)=p$.)

\item[III.] The only remaining case is where $A\in GL_2(\Z_p)$ has Jordan normal form $\tilde{A}$ $ \ =\left(\ba{cc}\beta & 0\\0&\ol{\beta} \ea\right)$, where $\beta,\ol{\beta}\in \Z_p(\beta)$, a quadratic field extension of $\Z_p$ (and $\beta$ and $\ol{\beta}$ are conjugate roots of an irreducible quadratic polynomial over $\Z_p$). For any $d$ such that $d\,|\,p^2-1$ but $d\nmid p-1$, there exists such an $A$ (and hence $\tilde{A}$) having order $d$.  

It follows that all cycles in $\Z_p(\beta)^2$, except the trivial one $\left(\ba{c} 0\\0 \ea\right)$, have length $d$.  Note that the orbit of $\tilde{A}$ on $\left(\ba{c} x\\y \ea\right)$,  for $x,y\in \Z_p(\beta)$, not both $0$, is given by   

$$\left(\ba{c} x\\y \ea\right), \left(\ba{c}\beta x  \\ \ol{\beta} y \ea\right),\left(\ba{c}\beta^2 x \\ \ol{\beta}^2 y \ea\right),\ldots, \left(\ba{c}\beta^{d-1} x\\ \ol{\beta}^{d-1} y \ea\right),$$

of length $d$, since $o^\times(\ol{\beta}) = o^\times(\beta) = d$. Hence, by Lemma 2.4, all nontrivial cycles of $A$ in $\Z_p^2$ also have length $d$.

\ee  

We are now ready to characterise all automorphisms of $\Z_p^2$.

{\bf 2.5. Theorem.\ } Let $|A|=p^2$ where $p$ is prime, and let $f:A\ra A$ be a bijection. Then $A$ can be turned into a group isomorphic to $\Z_p^2$, with $f\in \mb{Aut}(\Z_p^2)$, if and only if $f$ is the identity map, or $f$ has any one of the following cycle structures:
\be 
\item[(a)] $\left[\ba{cc} \frac{p^2-1}{d} & 1\\ d&1\ea\right]$ for some divisor $d$ of $p^2-1$;
\item[(b)] $\left[\ba{ccc} \frac{p-1}{d} & \frac{p-1}{d}& 1\\ pd&d&1\ea\right]$ for some divisor $d$ of $p-1$;
\item[(c)] $\left[\ba{cccc} \frac{(p-1)^2}{\mb{\footnotesize lcm}(d_1, d_2)} & \frac{p-1}{d_1}&\frac{p-1}{d_2}& 1\\ \mb{\footnotesize lcm}(d_1, d_2)&d_1&d_2&1\ea\right]$ for divisors $d_1$ and $d_2$ of $p-1$;
\ee

{\bf Proof:\ } From case III in the discussion preceding the theorem, an automorphism $f$ with cycle structure $\left[\ba{cc} 1 & 1\\ p^2-1&1\ea\right]$ exists.  By letting $d$ vary over all divisors of $p^2-1$, and by considering the corresponding automorphisms $f^d$ ($f$ composed with itself $d$ times), we obtain all the possible cycle structures given in (a). The cycle structures in (b) and (c) follow from cases II and I respectively. Also note that, if it happens that $d_1=d_2 = d$ in (c), where $d\,|\,p-1$, then the cycle structure in (a) is obtained. In particular, $d_1=d_2=1$ gives the identity map.\ \bb

{\bf 2.6. Example.\ } Let $p=7$. The divisors of $p^2-1=48$ that are not divisors of $p-1=6$, are given by $d\in \{4, 8,12,16,24,48\}$. For these divisors we obtain, from Theorem 2.5(a), the following corresponding cycle structures: 

$$\left[\ba{cc} 12 & 1\\ 4&1\ea\right], \ \left[\ba{cc} 6 & 1\\ 8&1\ea\right], \ \left[\ba{cc} 4 & 1\\ 12&1\ea\right], \ \left[\ba{cc} 3 & 1\\ 16&1\ea\right], \ \left[\ba{cc} 2 & 1\\ 24&1\ea\right], \ \left[\ba{cc} 1 & 1\\ 48&1\ea\right].$$

The divisors of $p-1=6$ are $d\in \{1,2,3,6\}$, so Theorem 2.5(b) gives the corresponding cycle structures

$$\left[\ba{ccc} 6 & 6 & 1\\ 7&1&1\ea\right]=\left[\ba{cc} 6 & 7\\ 7&1\ea\right], \ \left[\ba{ccc} 3 & 3 & 1\\ 14&2&1\ea\right], \ \left[\ba{ccc} 2 & 2 & 1\\ 21&3&1\ea\right], \ \left[\ba{ccc} 1 & 1 & 1\\ 42&6&1\ea\right].$$
 
Finally, for the remaining cases, we consider Theorem 2.5(c), where we take $d_1,d_2\in\{1,2,3,6\}$ and we may assume that $1\le d_1 \le d_2\le 6$. We obtain the cycles

$$\left[\ba{cccc} 36 & 6 & 6& 1\\ 1&1&1&1\ea\right] = \left[\ba{c} 49\\ 1\ea\right], \ \left[\ba{cccc} 18 & 6 & 3& 1\\ 2&1&2&1\ea\right]=\left[\ba{cc} 21 & 7\\ 2&1\ea\right], \ 
\left[\ba{cccc} 12 & 6 & 2& 1\\ 3&1&3&1\ea\right] = \left[\ba{cc} 14 & 7\\ 3&1\ea\right], $$
$$\left[\ba{cccc} 6 & 6 & 1& 1\\ 6&1&6&1\ea\right] = \left[\ba{cc} 7 & 7\\ 6&1\ea\right], \ 
\left[\ba{cccc} 18 & 3 & 3& 1\\ 2&2&2&1\ea\right] = \left[\ba{cc} 24 & 1\\ 2&1\ea\right], \ 
\left[\ba{cccc} 6 & 3 & 2& 1\\ 6&2&3&1\ea\right] = \left[\ba{cccc} 6 & 2&3&1\\ 6&3&2&1\ea\right], $$
$$ \left[\ba{cccc} 6 & 3 & 1& 1\\ 6&2&6&1\ea\right] = \left[\ba{ccc} 7 & 3&1\\ 6&2&1\ea\right], \ 
\left[\ba{cccc} 12 & 2 & 2& 1\\ 3&3&3&1\ea\right] = \left[\ba{cc} 16&1\\ 3&1\ea\right], \
\left[\ba{cccc} 6 & 2 & 1& 1\\ 6&3&6&1\ea\right] = \left[\ba{ccc} 7 & 2&1\\ 6&3&1\ea\right],  $$
$$ \left[\ba{cccc} 6 & 1 & 1& 1\\ 6&6&6&1\ea\right] = \left[\ba{cc} 8 & 1\\ 6&1\ea\right].$$

Consequently, if $|A|=49$, then $A$ can be made into a group isomorphic to $\Z_7^2$, with $f\in \mb{Aut}(\Z_7^2)$, if and only if $f$ has one of the $20$ cycle structures shown here.

One immediately raises the question of how the cycle structures of automorphisms on $\Z_{p^2}$ relate to the cycle structures of automorphisms on $\Z_p \times \Z_p$. It turns out that the former forms a subset of the latter.

{\bf 2.7. Theorem.\ } Let $|A|=p^2$, with $p$ prime, and let $f:A\ra A$ be a bijection.  Then $f$ has the auto-property if and only if $f$ has one of the cycle structures of an automorphism of $\Z_p^2$, given by Theorem 2.5.

{\bf Proof: \ } It suffices to show that every cycle structure that appears in Example 2.3(c), also appears in Theorem 2.5. For the parameters $\ell_i=\mb{ord}_{p^2}(k_i),\ 2\le i\le p^2-p$, and $\lambda = \mb{ord}_p(k_i)$, we have $\ell_i\mid \phi(p^2)$ and $\lambda\mid \phi(p)$ (see \cite[Theorem 2.14]{Nat}), and $\lambda\mid \ell_i$ (see Example 2.3(c)). So we have the two possibilities: 
\be 
\item $\ell_i=\lambda$. This gives the cycle structure $\left[\ba{cc} \frac{p^2-p}{\lambda}+\frac{p-1}{\lambda} & 1\\ \lambda &1\ea\right] = \left[\ba{cc} \frac{p^2-1}{\lambda} & 1\\ \lambda &1\ea\right]$, which agrees with the cycle structure of Theorem 2.5(c) with $d_1=d_2=\lambda$. 
\item $\ell_i = p\lambda$. This gives the cycle structure $\left[\ba{ccc} \frac{p^2-p}{p\lambda}&\frac{p-1}{\lambda} & 1\\ p\lambda & \lambda&1\ea\right] = \left[\ba{ccc} \frac{p-1}{\lambda}&\frac{p-1}{\lambda} & 1\\ p\lambda & \lambda&1\ea\right]$, which agrees with the cycle structure of Theorem 2.5(b) with $d=\lambda$.  \  \bb
 \ee

A natural question is whether this result holds in a more general setting, i.e., whether, for a given prime $p$ and an integer $n\ge2$, the cycle structures of the automorphisms of $\Z_p^n$ already contain all possible cycle structures of all abelian groups of order $p^n$. This is unfortunately not the case, as the next example shows: 

{\bf 2.8 Example.\ } The group $\mathbb{Z}_8$ has an automorphism of which the cycle structure is different from that of all automorphisms of $\mathbb{Z}_2^3$. 

{\bf Proof: \ } Consider $f:\mathbb{Z}_8\rightarrow \mathbb{Z}_8$ defined by $f(x)=-x$. The cycle structure of $f$ is $\left[\ba{cc} 2 & 3\\ 1 &2\ea\right]$. Assume that there is an automorphism $g$ of $\mathbb{Z}_2^3$ which has the same cycle structure. 

Consider the Jordan canonical form of $g$. Then $g$ has minimal polynomial $(x+1)^2$ since all elements of $\Z_2^3$ lie within cycles of length at most $2$. The characteristic polynomial of $g$ is therefore $(x+1)^3$, so that there are two blocks in Jordan form, one of size $2\times 2$ and one of size $1\times 1$. Consequently, the eigenspace related to the eigenvalue $-1 = 1$ must be of dimension at least $2$, implying that there will be at least four elements of $\Z_2^3$ in cycles of length $1$, a contradiction. \ \bb \\

{\bf 3. Groups isomorphic to $\Z^n$}

We will now investigate the cycle structures of the automorphisms of all groups of the form $\mathbb{Z}^n$. We must clearly still have the trivial cycle of length $1$, representing $0\mapsto 0$. From now on, we will refer to this cycle as the {\it zero cycle} of the map. Since $\mathbb{Z}$ is infinite there is the possibility of not only having (finite) cycles such as with the cases in Section 2, but also {\it chains}, i.e., distinct elements $\ldots, a_i, a_{i+1},a_{i+2},\ldots$ from $A$, such that $f(a_i) = a_{i+1}$ for all $i$. 

One of the major tools that we used to investigate the automorphisms of the finite groups was the fact that the elements of the general linear group were much more than just matrices over rings, but they were actually matrices over fields, which allowed us to use the Jordan normal form to form conjugacy classes which partitioned the general linear group. In the case of matrices over $\mathbb{Z}$ this cannot be done, as $\mathbb{Z}$ is not a field. But even though we have lost the Jordan normal forms, we still have that the automorphism group is isomorphic to $GL(\mathbb{Z},n)$. These are clearly all the $n\times n$ integer matrices with determinant equal to $\pm1 $ (\cite{Newman_IM}). 

For $1\le i \le n$,  $e_i$ is used to denote the element $(0,\ldots,1,\ldots,0) \in \Z^n$, with $1$ in the $i$-th coordinate and zeros elsewhere.   

{\bf 3.1 Proposition.\ } \label{NesZn} Suppose $f\in \mb{Aut}(\Z^n)$ for some positive integer $n$. Then the following conditions must hold:
\begin{enumerate}
\item If $f$ has a cycle, apart from the zero cycle, of any length $k$, then $f$ has infinitely many cycles of length $k$. 
\item If $f$ has a chain, it has infinitely many chains.  
\item If all $e_i,i\in\{1,2,\ldots,n\}$ are in cycles of $f$, then all elements of $\Z^n$ are in cycles, i.e., $f$ has no chains.  
\end{enumerate}

{\bf Proof:}\  
Let the matrix representation of $f$ be $M$, and represent the elements of the group $\Z^n$ as columns.
\begin{enumerate}
\item  Consider any non-zero cycle $T=(x,Mx,M^2x,\ldots,M^{k-1}x)$ of length $k$. Let $S_T$ be the (finite) set of all the absolute values of the non-zero components of the members of $T$. There exists a (non-zero) smallest element in $S_T$. Now, for any positive integer $m$, we see that $mT=(mx,M(mx),M^2(mx),\ldots,M^{k-1}(mx))=(mx,mMx,mM^2x,\ldots,mM^{k-1}x)$ is a cycle of length $k$, with $S_{mT}=mS_T$, from which it follows that the cycles $mT$ are disjoint for different $m\in\N$ as the minimum components are all distinct from one another. Consequently there are infinitely many cycles of length $k$.
  
\item The proof is roughly the same as above. The cycle $T=(x,Mx,\ldots,M^{k-1}x)$ is just replaced by the chain $T = (\ldots, M^{-2}x,M^{-1}x,x,Mx,M^2x,\ldots).$ Here,  the existence of the smallest (non-zero) element is guaranteed by the well-ordering principle on $\N$. 
\item Suppose all the $e_i$ are in cycles with the cycle containing $e_i$ of length $k_i$. Any $x\in\Z^n$ can be represented as $x=\sum_{i=1}^n\alpha_ie_i,\alpha_i\in\Z$. Denote the least common multiple of the set $\{k_j,j\in\{1,2,\ldots,n\}\}$ by $\ell$, and define $q_i=\frac{\ell}{k_i}$. Then
\[
M^\ell x  = \sum_{i=1}^n \alpha_iM^\ell e_i 
  = \sum_{i=1}^n \alpha_i e_i  = x,
\]

as $M^{\ell}e_i=e_i$ for all $i = 1,2,\ldots,n$. Hence $x$ lies in a cycle of length dividing $\ell$. \ \bb
\end{enumerate}

Proposition 3.1 tells us that if the structure of an automorphism consists of cycles only, then there is only a finite number of possible cycle lengths, as all cycles must be of length dividing the least common multiple of the lengths of the cycles of the $e_i$'s. However, it is still possible, in principle, for an infinite number of distinct cycle lengths to exist for an automorphism, but then some $e_i$ must lie in a chain. The following result shows that not even this is possible.\\

{\bf 3.2 Proposition.}\ 
The structure of any automorphism of $\Z^n$ possesses at most finitely many distinct cycle lengths. 

{\bf Proof:} \ 
Suppose that $f$ has infinitely many distinct cycle lengths. Let the matrix representation of $f$ be given  by the  $n\times n$ matrix $M$. 

If $n=1$, then $f(1)=1$ or $f(1)=-1$, as $\det(M)=\pm 1$. The former case is simply the identity mapping, and the latter has a cycle structure consisting of only cycles of length two, together with the zero-cycle. 

So assume that $n\geq 2$ for the remainder of the proof. Let, for $1\le i \le n$, 
 \[
x_i=\begin{bmatrix}
x_{1i}\\
x_{2i} \\
\vdots \\
x_{ni}
\end{bmatrix}
\]
be any $n$ distinct non-zero elements of $\Z^n$ occurring in cycles. For each $i$, let the cycle length of $x_i$ be $s_i$. Let $U = [x_1|x_2|\ldots|x_n]$, the $n\times n$ matrix with columns $x_i, i=1,\ldots,n$. From Proposition~3.1 it follows that at least one of the $e_i$'s belongs to a chain, say it is $e_1$. Also, let $U(r,c)$ denote the $(r,c)$-minor of $U$. We now consider 
\[y = \sum_{i=1}^n(-1)^{i-1}U(1,i)x_i.\]
For each $i\in\{2,\ldots,n\}$, the entry in the $i$-th row of $y$ is clearly the determinant of the matrix 
\begin{center}
$
\begin{bmatrix}
x_{i1} & x_{i2} & \ldots & x_{in} \\
x_{21} & x_{22} & \ldots & x_{2n} \\
\vdots & \vdots & \ddots & \vdots \\
x_{i1} & x_{i2} & \ldots & x_{in} \\
\vdots & \vdots & \ddots & \vdots \\
x_{n1} & x_{n2} & \ldots & x_{nn} 
\end{bmatrix}, $
\end{center}
which is zero, as the $i$-th row is identical to the first row. This means all the components of $y$, except perhaps the first, are equal to zero.

In the same way we see that the first component of $y$ is simply the determinant of $U$. By denoting the least common multiple of $\{s_i, i \in\{1,2,\ldots,n\}\}$ by $\ell$, it is clear that
\[
M^\ell y  = \sum_{i=1}^n (-1)^{i-1}U(1,i)M^\ell x_i \\
 = \sum_{i=1}^n (-1)^{i-1}U(1,i)x_i \\
 = y,
\]
which means that $y$ belongs to a cycle. However, the element $e_1$ lies in a chain, implying that all non-zero elements with only their first components non-zero, lie in a chain. Consequently, $y$ must be $0$, from which it follows that $\det(U)=0$. This means that the columns of $U$ are  linearly dependent, and for fixed $x_1,x_2,\ldots,x_{n-1}$, any other element $z$ that lies in some cycle, can be expressed as $z = \sum_{i=1}^{n-1} \gamma_ix_i$, with $\gamma_i\in\Q$. Hence $z$ must belong to a cycle having a length dividing the least common multiple of the set $\{s_1,s_2,\ldots,s_{n-1}\}$. Since this holds for all $z$ occurring in cycles, we see that the structure of any automorphism has only finitely many distinct cycle lengths. \ \bb

{\bf 3.3 Example.} \ The structure of the automorphism on $\Z^2$ represented by $M=\left[\ba{cc}1&1\\1&0\ea\right]$ does not have any non-zero cycles, hence it consists only of the zero cycle, and infinitely many chains. 

{\bf Proof:}\ First we notice that for any positive integer $n$, $M^n= \left[\ba{cc}F_{n+1}&F_n\\F_n&F_{n-1}\ea\right]$ with $F_n$ the $n$-th number in the Fibonacci sequence. Suppose that the structure of $M$ contains a cycle of length $n\in\N$. Then there exist $a,b\in\Z$ such that
\begin{align*}
aF_{n+1} + bF_n &= a \\
aF_n+bF_{n-1} &= b,
\end{align*} 
which can be written as 
\begin{align*}
(F_{n+1}-1)a + F_nb &= 0 \\
F_na+(F_{n-1}-1)b &= 0.
\end{align*}
The determinant of this system is $(F_{n+1}-1)(F_{n-1}-1)-F_n^2$, which reduces to $(F_{n+1}^2-F_{n+1}F_n-F_n^2) + 1 - (F_{n+1}+F_{n-1})$. Using the identity $F_{n+1}^2-F_{n+1}F_n-F_n^2= (-1)^n$, we see that the system has a non-zero determinant, and conclude that $\begin{bmatrix} a \\ b \end{bmatrix}=\begin{bmatrix} 0 \\ 0 \end{bmatrix}$ is the only solution.\ \bb

{\bf 3.4 Definition.} \ 
Suppose the structure of an automorphism of $\Z^n$ contains a cycle of length $k$. Then this cycle is called a {\it primitive cycle} of the structure if for any proper divisor $d$ of $k$, there are no non-zero cycles of length $d$ in this structure. In this case, we call $k$ a {\it primitive cycle length} of the structure of the automorphism.

We shall now investigate whether for any natural number $k$, there exists an automorphism for which all the non-zero cycles are of length $k$. In order to do so, we shall first take some inspiration on the construction of cycles from larger ones. Suppose, for example, the automorphism $M$ of $\Z^n$ has a cycle $(x,Mx,M^2x,\ldots,M^5x)$ of length $6$.

The cycle generated by 
$x+M^2x+M^4x$ is $(x+M^2x+M^4x,Mx+M^3x+M^5x)$, hence has length $1$ or $2$. Similarly, the cycle $(x+M^3x,Mx+M^4x,M^2x+M^5x)$, generated by $x+M^3x$, must have a length that divides $3$.

Note that these cycle lengths are not necessarily of lengths $2$ and $3$ respectively.  They could also be of length $1$. At first it seems that this could severely restrict the possibilities on the numbers which could be primitive cycle lengths. However, surprisingly, this result does not eventually restrict the numbers which are primitive cycle lengths, but rather tells us how to construct automorphisms with exactly those primitive cycle lengths. For our $6$-cycle case, for example, if we can somehow find an invertible integer matrix $M$ such that $I+M^2+M^4=I+M^3=0$, then the constructed elements which could have cycle lengths of $2$ and $3$ will actually turn out to be the zero element, and the cycle reduces to the zero-cycle. 

{\bf 3.5 Example.}\ The automorphism on $\Z^2$ with matrix representation $M=\left[\ba{cc}0&1\\-1&1\ea\right]$ has all of its non-zero cycle lengths equal to $6$. 

{\bf Proof: \ }
We have that $I+M^2+M^4=0$ and $I+M^3=0$, and also that $M^6=I$, which means that all cycles are of length dividing $6$. Consider an arbitrary $\begin{bmatrix}x \\ y\end{bmatrix}\in\Z^2$. This  gives the cycle
\begin{center}
$\left( \begin{bmatrix}x \\ y\end{bmatrix}, \begin{bmatrix}y \\ y-x\end{bmatrix}, \begin{bmatrix}y-x \\ -x\end{bmatrix}, \begin{bmatrix}-x \\ -y\end{bmatrix}, \begin{bmatrix}-y \\ x-y\end{bmatrix}, \begin{bmatrix}x-y \\ x\end{bmatrix} \right)$.
\end{center} 
This is an explicit example of an automorphism of which the structure consists of the zero cycle, no chains and all non-zero cycles of length $6$, implying that $6$ is a primitive length with respect to this structure. Note that if any of these cycles were to collapse into a cycle of length less than $6$, then we must have that  
\begin{center}
$\begin{bmatrix}x \\ y\end{bmatrix}\in 
\left\{  \begin{bmatrix}y \\ y-x\end{bmatrix}, \begin{bmatrix}-x \\ -y\end{bmatrix}, \begin{bmatrix}y-x \\ -x\end{bmatrix}\right\}.$
\end{center} 
All these possibilities lead to the zero cycle. \ \bb

This example also paves the way towards establishing a technique that will allow us, for any positive integer $k$, the construction of an automorphism on some $\Z^n$ of which the structure has all of its non-zero cycles of length $k$. The next theorem is the first step towards this goal:

{\bf 3.6 Theorem.\ } For any $n>1$, let $n=\prod_{i=1}^k p_i^{\alpha_i}$ be the prime factorization of $n$, where we assume that $p_1>p_2>\cdots >p_k$. Define, for each $i\in\{1,\ldots, k\}$ the polynomial $Q_i$ by 
\[Q_i(\lambda)=\sum_{j=0}^{p_i-1}\lambda^{\frac{n\cdot j}{p_i}}.\]

Then the $n$-th cyclotomic polynomial $\Phi_n$ divides $Q_i$ for all $i\in\{1,2,\ldots,k\}$. Moreover, $\Phi_n$ is the only non-constant polynomial that divides all the $Q_i$.

{\bf Proof: \ }
First we notice that $\lambda^n-1=(\lambda^{\frac{n}{p_i}}-1)Q_i$ for any $i\in\{1,2,\ldots,k\}$. Let $\zeta$ be a primitive $n$-th root of unity. From $(\zeta^{\frac{n}{p_i}}-1)Q_i(\zeta)=0$ and  $\zeta^{\frac{n}{p_i}}-1\neq 0$ it follows that $Q_i(\zeta)=0$. An immediate consequence  is that $\lambda-\zeta$ is a factor of $Q_i$ for all primitive roots $\zeta$ of unity, so the $n$-th cyclotomic polynomial $\Phi_n$ divides all of the $Q_i$. 

Now suppose that there is another non-constant polynomial $R$ which is a factor of all the $Q_i$'s but with a root $\eta$ which is not a primitive $n$-th root of unity. As the roots of $R$ must all be $n$-th roots of unity, it follows that $\eta=\zeta^m$ for some $m\in\{1,2,\ldots,n\}$ and such that $\mbox{gcd}(m,n)\neq 1$. However, then there exists an $i$ such that $\eta^\frac{n}{p_i}-1=0$, and as $Q_i(\eta)=0$, it follows that $\eta$ is a root of $\lambda^n-1$ of multiplicity at least two. This is a contradiction, as all roots of $\lambda^n-1$ have multiplicity $1$. \ \bb

%

Combining Theorem 3.6 and The Cayley-Hamilton Theorem, it is clear that 
if we can find an $n\times n$ matrix $M$ with characteristic polynomial $\Phi_n$, then $M$ is a root of $\Phi_n$, and thus of all the $Q_i$'s.

We now have:

{\bf 3.7 Proposition.\ } Let $n\in \N$. Then there exists an automorphism $f_n:\mathbb{Z}^m\rightarrow\mathbb{Z}^m$, for some positive integer $m$, such that the structure of $f_n$ consists of only the zero cycle and infinitely many cycles of length $n$. 

{\bf Proof:\ }
Theorem 3.6 shows that the $n$-th cyclotomic  polynomial is the (non-constant) greatest common divisor of the $Q_i$'s.  Let $C_{\Phi_n}$ be its  companion matrix (so that $C_{\Phi_n}$ has characteristic polynomial $\Phi_n$). Since the constant term of $\Phi_n$ is either $1$ or $-1$, we have that $\det C_{\Phi_n}=\pm1$. Hence $C_{\Phi_n}$ is invertible, making it the matrix representation of an automorphism. $C_{\Phi_n}$ is a root of $\Phi_n$, and since $\Phi_n$ divides all the $Q_i$'s, we have that $C_{\Phi_n}$ is a root of all the $Q_i$'s.   

Since all the $Q_i$'s divide $\lambda^n-1$,  all cycles associated with $C_{\Phi_n}$ have lengths dividing $n$. Any cycle length $d$ properly dividing $n$, would have to divide $\frac{n}{p_i}$ for some $p_i$. By letting $\begin{bmatrix} x \\ y \end{bmatrix}$ be a non-zero element in any cycle of length $d$, we note that 
\[
Q_i(C_{\Phi_n})\begin{bmatrix}x \\y \end{bmatrix}  =\sum_{j=0}^{p_i-1}(C_{\Phi_n})^{\frac{nj}{p_i}}\begin{bmatrix}x \\y \end{bmatrix} 
  = \sum_{j=0}^{p_i-1} \begin{bmatrix}x \\y \end{bmatrix} 
  = p_i\begin{bmatrix}x \\y \end{bmatrix} 
  = \begin{bmatrix}p_i x \\p_iy \end{bmatrix},
\]
which clearly cannot hold, since $Q_i(C_{\Phi_n})=0$. The automorphism on $\Z^m$, where $m=\phi(n)$, of which $C_{\Phi_n}$ is the matrix representation consequently has a structure consisting of the zero-cycle, no chains, and only cycles of length $n$. Note, we cannot use Theorem 3.6 if $n=1$, but, of course, an automorphism with all its (non-zero) cycles of length $1$ does exist -- simply take the identity map on the group $\Z^m$, for any $m\ge1$.\ \bb

For each $n\in \N$, we shall call the automorphism described in Proposition 3.7 a {\it pure $n$-cyclic automorphism} and denote its matrix representation by $P_n$. 

We now proceed to investigate automorphisms on $\Z^n$ with cycles of different lengths.

{\bf 3.8 Theorem.\ }
Suppose the structure of an automorphism on $\Z^n$ has non-zero cycles of lengths $\alpha$ and $\beta$. Then the structure also has a cycle of length $[\alpha,\beta]$ (the least common multiple of $\alpha$ and $\beta$). 

{\bf Proof:\ }
Let $M$ be the matrix representation of the automorphism. Suppose $x$ lies in a cycle of length $\alpha$ and $y$ in a cycle of length $\beta$. It is clear that for each positive integer $k$, 
\[
M^{[\alpha,\beta]}(x+ky)  = M^{[\alpha,\beta]}x + kM^{[\alpha,\beta]}y 
 = x+ky,
\]
as $\alpha|[\alpha,\beta]$ and $\beta|[\alpha,\beta]$. Denote the cycle length of $x+ky$ by $\gamma_k$ for all $k\in\N$. Clearly, $\gamma_k|[\alpha,\beta]$, so there exist distinct $k,j\in\N$ with $\gamma_k=\gamma_j$. Denote this common value by $\gamma$. Consider the two cycles $(x+ky,M(x+ky),\ldots, M^{\gamma-1}(x+ky))$ and $(x+jy,M(x+jy),\ldots, M^{\gamma-1}(x+jy))$. Since matrix multiplication is distributive over matrix summation we can subtract these two cycles term by term to obtain a new cycle $((j-k)y,M(j-k)y,\ldots,M^{\gamma-1}(j-k)y)$. Note though, that the cycle length of $(j-k)y$ need not be $\gamma$. It is possible that the newly formed cycle actually fully traverses the cycle containing $(j-k)y$ several times. However, the cycle length of $(j-k)y$ must divide $\gamma$. Since $k\neq j$ and $M(j-k)y=(j-k)My$, it is clear that $(j-k)y$ must be in a cycle of the same length as $y$, and hence $\beta|\gamma$.

Let $\gamma=q\alpha+r,\  0\leq r<\alpha$. Then
\[
x+ky = M^\gamma(x+ky) =M^rx+kM^\gamma y 
 = M^rx + ky. \\
\]  
So $M^rx=x$, but since the cycle containing $x$ is of length $\alpha$, it follows that $r=0$, and so $\alpha\mid\gamma$. It follows that $[\alpha,\beta]|\gamma$, and we conclude that $\gamma=[\alpha,\beta]$. \ \bb

We can now give a complete structural characterization of all functions having the auto-property with underlying group $\Z^n$. 

{\bf 3.9 Theorem.\ } Let $A$ be a countably infinite set.  
A bijective function $f:A\rightarrow A$ possesses the auto-property with underlying group structure $(\Z^n,+)$ (for some $n\ge1$) if and only if the structure of $f$ satisfies all of the following:

\begin{enumerate}
\item[(1)] It contains at least one cycle of length $1$.
\item[(2)] The number of distinct cycle lengths is finite. 
\item[(3)] If it contains a non-zero cycle then it contains infinitely many cycles of this length.
\item[(4)] If it contains a chain, it contains infinitely many chains.
\item[(5)] If it contains non-zero cycles of length $\alpha$ and $\beta$, then it contains a cycle of length $[\alpha,\beta]$.
\end{enumerate}

{\bf Proof:\ }
Propositions 3.1, 3.2 and Theorem 3.8 show that the conditions listed above are necessary. 

Let $f$ be a function on a countably infinite set satisfying all the conditions listed in the theorem. We show that $f$ has the auto-property by constructing an invertible integer matrix representing $f$. Condition (2) allows the construction of a finite set $\mathcal{L}=\{n_1,n_2,\ldots,n_s\}$ consisting of the distinct cycle lengths occurring in the structure of $f$. For each $n_i\in\mathcal{L}$, Proposition 3.7 shows the existence of a pure $n_i$-cyclic automorphism. If $f$ has no chains, construct the integer matrix 

\[M=\begin{bmatrix}
P_{n_1} & 0 & 0 & \cdots & 0 \\
0 & P_{n_2} & 0 & \cdots & 0 \\
0 & 0 & P_{n_3}& \cdots & 0 \\
\vdots & \vdots & \vdots & \ddots & \vdots \\
0 & 0 & 0 & \cdots & P_{n_s}
\end{bmatrix}\]  

which is a block diagonal matrix obtained by placing the matrices $P_{n_i}$, as defined in Proposition 3.7 (as blocks) along the diagonal of $M$ and all other entries equal to $0$. If $f$ has chains, simply append the matrix $\left[\ba{cc}1&1\\1&0\ea\right]$ along the diagonal of $M$, say at the bottom on the right. 

Since all of the $P_{n_i}$'s are along the diagonal, it follows that $\det(M) = \det(P_{n_1})\det(P_{n_2})\cdots \det(P_{n_s})$ is either $1$ or $-1$, as all the $P_{n_i}$'s are invertible. This shows that $M$ is invertible, and represents an automorphism $f_M:\Z^m\rightarrow\Z^m$ for some positive integer $m$. Let the number of rows of $P_{n_i}$ be denoted by $x_i$. For $n_i$, it is clear that the cycle containing the element $e_{x_1+\cdots+x_{i-1}+1}$ is of length $n_i$ in the structure of $f_M$, as the cycle of $e_1$ is of length $n_i$ in the structure of the pure ${n_i}$-cycle represented by the matrix $P_{n_i}$. The structure of $f_M$ thus contains cycles of length $n_i$ for each $n_i\in\mathcal{L}$. If $f$ contains a chain, the last matrix embedded in the diagonal of $M$ is $\left[\ba{cc}1&1\\1&0\ea\right]$. Example 3.3 then shows that $e_{x_1+\cdots+x_s+1}$  lies in a chain. It is now clear that a non-zero cycle of length $n_i$ (or a chain) occurs in the structure of $f$ only if one also occurs in the structure of $f_M$. 

Given any element $z\in\Z^m$, written as a column, we can decompose $z$ as the sum $z=z_1+z_2+\cdots+z_s+\hat{z}$ with each $z_n$ being a column of length $m$, with its $j$-th entry equal to that of $z$, for all $j\in\left\{\left( \sum_{i=1}^{n-1}x_i\right) +1,\ldots,\left( \sum_{i=1}^{n-1}x_i\right)+x_n \right\}$, and zeros elsewhere. If $f$ has chains, $\hat{z}$ is a column of length $m$, with the first $m-2$ entries equal to $0$, and the last two entries equal to the corresponding entries of $z$; otherwise put $\hat{z}$ equal to the zero column of length $m$,  i.e., all its entries are equal to $0$. We will refer to $z_i$ as the $n_i$-cycle component of $z$, and to $\hat{z}$ as the chain component of $z$. Since $M^{\ell}z=\sum_{j=1}^s M^{\ell}z_j+ M^{\ell}\hat{z}$ for all $\ell\ge1$, and each $z_i$ is in a cycle of length dividing $n_i$ \ ($1\le i \le s$), it is clear that $z$ is in a chain if and only if $\hat{z}$ is in a chain, which is the case for exactly all non-zero $\hat{z}$. Consequently, if $f_M$ has a chain, then $f$ must also have had one (since otherwise $\hat{z}=0$ for all $z\in\Z^m$). Now take any $z$ in a non-zero cycle of $f_M$. As discussed above, $\hat{z}$ must be the zero column. However, since $M$ acts on $z_k$ in the same way as the pure $n_k$-cycle would on a column consisting of the $\left( \left( \sum_{i=1}^{k-1}x_i\right) +1\right)$-th up to $\left( \left( \sum_{i=1}^{k-1}x_i\right) +x_k\right)$-th entries of $z_k$, it follows that the cycle of $z_k$ is either the zero-cycle, or of length $n_k$. Since the $z_i$'s are linearly independent, the cycle length of $z$ is equal to the least common multiple of the $n_i$'s for which the corresponding $z_i$'s are not zero columns. It now follows that any non-zero cycle of $f_M$ has length the least common multiple of  $n_{\sigma(1)},n_{\sigma(2)},\ldots,n_{\sigma(k)}$ for some permutation $\sigma$ of $(1,2,\ldots,s)$, with $k\leq s $, and (by condition $5$) of the same length as some cycle of $f$. Consequently, a cycle of length $n$ (or a chain) occurs in the structure of $f_M$ only if one also occurs in that of $f$.  We now have that the structures of $f$ and $f_M$ have cycles of the same distinct lengths (as well as chains) if and only if the other one has, and by conditions (1), (3) and (4), infinitely many of them, apart from the  zero-cycle. It follows that $f$ has the auto-property.  \ \bb

{\bf 3.11 Example.} \ Suppose we want to construct a matrix which represents an automorphism with chains, and cycles of lengths $6$ and $15$. 

As there are cycles of length $6$ and $15$, there must be a cycle of length $30$. We proceed to find $P_6,P_{15}$ and $P_{30}$.

 {\bf Pure $\mathbf{6}$-cycle}: $Q_1(\lambda)=1+\lambda^2+\lambda^4$ and $Q_2(\lambda)=1+\lambda^3$. The $\gcd$ of the $Q_i$'s is  $\Phi_6(\lambda)=1-\lambda+\lambda^2$. The companion matrix of this polynomial is:
\[P_6=\begin{bmatrix}
\phantom{-}0 & 1 \\
-1 & 1
\end{bmatrix}.\]
{\bf Pure $\mathbf{15}$-cycle}: $Q_1(\lambda)=1+\lambda^3+\lambda^6+\lambda^9+\lambda^{12}$ and $Q_2(\lambda)=1+\lambda^5+\lambda^{10}$. The $\gcd$ of the $Q_i$'s is $\Phi_{15}(\lambda)=1-\lambda+\lambda^3-\lambda^4+\lambda^5-\lambda^7+\lambda^8$. The companion matrix of this polynomial is:
\[P_{15}=\begin{bmatrix}
\phantom{-}0 & \phantom{-}1 & \phantom{-}0 & \phantom{-}0 & \phantom{-}0 & \phantom{-}0 & \phantom{-}0 & \phantom{-}0\\
\phantom{-}0 & \phantom{-}0 & \phantom{-}1 & \phantom{-}0 & \phantom{-}0 & \phantom{-}0 & \phantom{-}0 & \phantom{-}0\\
\phantom{-}0 & \phantom{-}0 & \phantom{-}0 & \phantom{-}1 & \phantom{-}0 & \phantom{-}0 & \phantom{-}0 & \phantom{-}0\\
\phantom{-}0 & \phantom{-}0 & \phantom{-}0 & \phantom{-}0 & \phantom{-}1 & \phantom{-}0 & \phantom{-}0 & \phantom{-}0\\
\phantom{-}0 & \phantom{-}0 & \phantom{-}0 & \phantom{-}0 & \phantom{-}0 & \phantom{-}1 & \phantom{-}0 & \phantom{-}0\\
\phantom{-}0 & \phantom{-}0 & \phantom{-}0 & \phantom{-}0 & \phantom{-}0 & \phantom{-}0 & \phantom{-}1 & \phantom{-}0\\
\phantom{-}0 & \phantom{-}0 & \phantom{-}0 & \phantom{-}0 & \phantom{-}0 & \phantom{-}0 & \phantom{-}0 & \phantom{-}1\\
-1 & \phantom{-}1 & \phantom{-}0 & -1 & \phantom{-}1 & -1 & \phantom{-}0 & \phantom{-}1 
\end{bmatrix}.\]  
{\bf Pure $\mathbf{30}$-cycle}: $Q_1(\lambda)=1+\lambda^6+\lambda^{12}+\lambda^{18}+\lambda^{24}$, $Q_2(\lambda)=1+\lambda^{10}+\lambda^{20}$ and $Q_3(\lambda)=1+\lambda^{15}$. The $\gcd$ of the $Q_i$'s is $\Phi_{30}(\lambda)=1+\lambda-\lambda^3-\lambda^4-\lambda^5+\lambda^7+\lambda^8$. The companion matrix of this polynomial is:
\[P_{30}=\begin{bmatrix}
\phantom{-}0 & \phantom{-}1 & \phantom{-}0 & \phantom{-}0 & \phantom{-}0 & \phantom{-}0 & \phantom{-}0 & \phantom{-}0\\
\phantom{-}0 & \phantom{-}0 & \phantom{-}1 & \phantom{-}0 & \phantom{-}0 & \phantom{-}0 & \phantom{-}0 & \phantom{-}0\\
\phantom{-}0 & \phantom{-}0 & \phantom{-}0 & \phantom{-}1 & \phantom{-}0 & \phantom{-}0 & \phantom{-}0 & \phantom{-}0\\
\phantom{-}0 & \phantom{-}0 & \phantom{-}0 & \phantom{-}0 & \phantom{-}1 & \phantom{-}0 & \phantom{-}0 & \phantom{-}0\\
\phantom{-}0 & \phantom{-}0 & \phantom{-}0 & \phantom{-}0 & \phantom{-}0 & \phantom{-}1 & \phantom{-}0 & \phantom{-}0\\
\phantom{-}0 & \phantom{-}0 & \phantom{-}0 & \phantom{-}0 & \phantom{-}0 & \phantom{-}0 & \phantom{-}1 & \phantom{-}0\\
\phantom{-}0 & \phantom{-}0 & \phantom{-}0 & \phantom{-}0 & \phantom{-}0 & \phantom{-}0 & \phantom{-}0 & \phantom{-}1\\
-1 & -1 & \phantom{-}0 & \phantom{-}1 & \phantom{-}1 & \phantom{-}1 & \phantom{-}0 & -1  
\end{bmatrix}.\]  

The matrix which represents the desired automorphism (on $\Z^{20}$) is

\begin{center}
\resizebox{.9\hsize}{!}{
$\begin{bmatrix}
 \phantom{-}0 & \phantom{-}1 & \phantom{-}0 & \phantom{-}0 & \phantom{-}0 & \phantom{-}0 & \phantom{-}0 & \phantom{-}0 & \phantom{-}0 & \phantom{-}0 & \phantom{-}0 & \phantom{-}0 & \phantom{-}0 & \phantom{-}0 & \phantom{-}0 & \phantom{-}0 & \phantom{-}0 & \phantom{-}0 & \phantom{-}0 & \phantom{-}0 \\ 
 -1 & \phantom{-}1 &  \phantom{-}0 &  \phantom{-}0 &  \phantom{-}0 &  \phantom{-}0 &  \phantom{-}0 &  \phantom{-}0 &  \phantom{-}0 & \phantom{-}0 &  \phantom{-}0 &  \phantom{-}0 &  \phantom{-}0 &  \phantom{-}0 &  \phantom{-}0 &  \phantom{-}0 &  \phantom{-}0 &  \phantom{-}0 &  \phantom{-}0 &  \phantom{-}0 \\
  \phantom{-}0 &  \phantom{-}0 &  \phantom{-}0 & \phantom{-}1 &  \phantom{-}0 &  \phantom{-}0 &  \phantom{-}0 &  \phantom{-}0 &  \phantom{-}0 &  \phantom{-}0 &  \phantom{-}0 &  \phantom{-}0 &  \phantom{-}0 &  \phantom{-}0 &  \phantom{-}0 &  \phantom{-}0 &  \phantom{-}0 &  \phantom{-}0 &  \phantom{-}0 &  \phantom{-}0\\
  \phantom{-}0 &  \phantom{-}0 &  \phantom{-}0 &  \phantom{-}0 & \phantom{-}1 &  \phantom{-}0 &  \phantom{-}0 &  \phantom{-}0 &  \phantom{-}0 &  \phantom{-}0 &  \phantom{-}0 &  \phantom{-}0 &  \phantom{-}0 &  \phantom{-}0 &  \phantom{-}0 &  \phantom{-}0 &  \phantom{-}0 &  \phantom{-}0 &  \phantom{-}0 &  \phantom{-}0\\
  \phantom{-}0 &  \phantom{-}0 &  \phantom{-}0 &  \phantom{-}0 &  \phantom{-}0 & \phantom{-}1 &  \phantom{-}0 &  \phantom{-}0 &  \phantom{-}0 &  \phantom{-}0 &  \phantom{-}0 &  \phantom{-}0 &  \phantom{-}0 &  \phantom{-}0 &  \phantom{-}0 &  \phantom{-}0 &  \phantom{-}0 &  \phantom{-}0 &  \phantom{-}0 &  \phantom{-}0\\
  \phantom{-}0 &  \phantom{-}0 &  \phantom{-}0 &  \phantom{-}0 &  \phantom{-}0 &  \phantom{-}0 & \phantom{-}1 &  \phantom{-}0 &  \phantom{-}0 &  \phantom{-}0 &  \phantom{-}0 &  \phantom{-}0 &  \phantom{-}0 &  \phantom{-}0 &  \phantom{-}0 &  \phantom{-}0 &  \phantom{-}0 &  \phantom{-}0 &  \phantom{-}0 &  \phantom{-}0\\
  \phantom{-}0 &  \phantom{-}0 &  \phantom{-}0 &  \phantom{-}0 &  \phantom{-}0 &  \phantom{-}0 &  \phantom{-}0 & \phantom{-}1 &  \phantom{-}0 &  \phantom{-}0 &  \phantom{-}0 &  \phantom{-}0 &  \phantom{-}0 &  \phantom{-}0 &  \phantom{-}0 &  \phantom{-}0 &  \phantom{-}0 &  \phantom{-}0 &  \phantom{-}0 &  \phantom{-}0\\
  \phantom{-}0 &  \phantom{-}0 &  \phantom{-}0 &  \phantom{-}0 &  \phantom{-}0 &  \phantom{-}0 &  \phantom{-}0 &  \phantom{-}0 & \phantom{-}1 &  \phantom{-}0 &  \phantom{-}0 &  \phantom{-}0 &  \phantom{-}0 &  \phantom{-}0 &  \phantom{-}0 &  \phantom{-}0 &  \phantom{-}0 &  \phantom{-}0 &  \phantom{-}0 &  \phantom{-}0\\
  \phantom{-}0 &  \phantom{-}0 &  \phantom{-}0 &  \phantom{-}0 &  \phantom{-}0 &  \phantom{-}0 &  \phantom{-}0 &  \phantom{-}0 &  \phantom{-}0 & \phantom{-}1 &  \phantom{-}0 &  \phantom{-}0 &  \phantom{-}0 &  \phantom{-}0 &  \phantom{-}0 &  \phantom{-}0 &  \phantom{-}0 &  \phantom{-}0 &  \phantom{-}0 &  \phantom{-}0\\
  \phantom{-}0 &  \phantom{-}0 & -1 & \phantom{-}1 &  \phantom{-}0 & -1 & \phantom{-}1 & -1 &  \phantom{-}0 & \phantom{-}1 &  \phantom{-}0 &  \phantom{-}0 &  \phantom{-}0 &  \phantom{-}0 &  \phantom{-}0 &  \phantom{-}0 &  \phantom{-}0 &  \phantom{-}0 &  \phantom{-}0 &  \phantom{-}0 \\ 
  \phantom{-}0 &  \phantom{-}0 &  \phantom{-}0 &  \phantom{-}0 &  \phantom{-}0 &  \phantom{-}0 &  \phantom{-}0 &  \phantom{-}0 &  \phantom{-}0&  \phantom{-}0 &  \phantom{-}0 & \phantom{-}1 &  \phantom{-}0 &  \phantom{-}0 &  \phantom{-}0 &  \phantom{-}0 &  \phantom{-}0 &  \phantom{-}0 &  \phantom{-}0 &  \phantom{-}0\\
  \phantom{-}0 &  \phantom{-}0 &  \phantom{-}0 &  \phantom{-}0 &  \phantom{-}0 &  \phantom{-}0 &  \phantom{-}0 &  \phantom{-}0 &  \phantom{-}0&  \phantom{-}0 &  \phantom{-}0 &  \phantom{-}0 & \phantom{-}1 &  \phantom{-}0 &  \phantom{-}0 &  \phantom{-}0 &  \phantom{-}0 &  \phantom{-}0 &  \phantom{-}0 &  \phantom{-}0\\
  \phantom{-}0 &  \phantom{-}0 &  \phantom{-}0 &  \phantom{-}0 &  \phantom{-}0 &  \phantom{-}0 &  \phantom{-}0 &  \phantom{-}0 &  \phantom{-}0&  \phantom{-}0 &  \phantom{-}0 &  \phantom{-}0 &  \phantom{-}0 & \phantom{-}1 &  \phantom{-}0 &  \phantom{-}0 &  \phantom{-}0 &  \phantom{-}0 &  \phantom{-}0 &  \phantom{-}0\\
  \phantom{-}0 &  \phantom{-}0 &  \phantom{-}0 &  \phantom{-}0 &  \phantom{-}0 &  \phantom{-}0 &  \phantom{-}0 &  \phantom{-}0 &  \phantom{-}0&  \phantom{-}0 &  \phantom{-}0 &  \phantom{-}0 &  \phantom{-}0 &  \phantom{-}0 & \phantom{-}1 &  \phantom{-}0 &  \phantom{-}0 &  \phantom{-}0 &  \phantom{-}0 &  \phantom{-}0\\
  \phantom{-}0 &  \phantom{-}0 &  \phantom{-}0 &  \phantom{-}0 &  \phantom{-}0 &  \phantom{-}0 &  \phantom{-}0 &  \phantom{-}0 &  \phantom{-}0&  \phantom{-}0 &  \phantom{-}0 &  \phantom{-}0 &  \phantom{-}0 &  \phantom{-}0 &  \phantom{-}0 & \phantom{-}1 &  \phantom{-}0 &  \phantom{-}0 &  \phantom{-}0 &  \phantom{-}0\\
  \phantom{-}0 &  \phantom{-}0 &  \phantom{-}0 &  \phantom{-}0 &  \phantom{-}0 &  \phantom{-}0 &  \phantom{-}0 &  \phantom{-}0 &  \phantom{-}0&  \phantom{-}0 &  \phantom{-}0 &  \phantom{-}0 &  \phantom{-}0 &  \phantom{-}0 &  \phantom{-}0 &  \phantom{-}0 & \phantom{-}1 &  \phantom{-}0 &  \phantom{-}0 &  \phantom{-}0\\
  \phantom{-}0 &  \phantom{-}0 &  \phantom{-}0 &  \phantom{-}0 &  \phantom{-}0 &  \phantom{-}0 &  \phantom{-}0 &  \phantom{-}0 &  \phantom{-}0&  \phantom{-}0 &  \phantom{-}0 &  \phantom{-}0 &  \phantom{-}0 &  \phantom{-}0 &  \phantom{-}0 &  \phantom{-}0 &  \phantom{-}0 & \phantom{-}1 &  \phantom{-}0 &  \phantom{-}0\\
  \phantom{-}0 &  \phantom{-}0 &  \phantom{-}0 &  \phantom{-}0 &  \phantom{-}0 &  \phantom{-}0 &  \phantom{-}0 &  \phantom{-}0 &  \phantom{-}0&  \phantom{-}0 & -1 & -1 &  \phantom{-}0 & \phantom{-}1 & \phantom{-}1 & \phantom{-}1 &  \phantom{-}0 & -1 &  \phantom{-}0 &  \phantom{-}0  \\
  \phantom{-}0 &  \phantom{-}0 &  \phantom{-}0 &  \phantom{-}0 &  \phantom{-}0 &  \phantom{-}0 &  \phantom{-}0 &  \phantom{-}0 &  \phantom{-}0 &  \phantom{-}0 &  \phantom{-}0 &  \phantom{-}0 &  \phantom{-}0 &  \phantom{-}0 &  \phantom{-}0 &  \phantom{-}0 &  \phantom{-}0 &  \phantom{-}0 & \phantom{-}1 & \phantom{-}1 \\
  \phantom{-}0 &  \phantom{-}0 &  \phantom{-}0 &  \phantom{-}0 &  \phantom{-}0 &  \phantom{-}0 &  \phantom{-}0 &  \phantom{-}0 &  \phantom{-}0 &  \phantom{-}0 &  \phantom{-}0 &  \phantom{-}0 &  \phantom{-}0 &  \phantom{-}0 &  \phantom{-}0 &  \phantom{-}0 &  \phantom{-}0 &  \phantom{-}0 & \phantom{-}1 &  \phantom{-}0 \\
\end{bmatrix}.$}
\end{center} 

{\bf Acknowledgement.}\ The first, second and fourth authors would like to thank the National Research Foundation of South Africa for financial assistance. All the authors would also like to thank the referees for valuable comments and suggestions.



\begin{thebibliography}{99}
\bibitem{glg} http://groupprops.subwiki.org/wiki/Element\_structure\_of\_general\_linear\_group\_of\_degree\_two\_\lb over\_a\_finite\_field
\bibitem{BouFonKeYeh} P. Bouchard, Y. Fong, W.-F. Ke and Y.-N. Yeh, Counting $f$ such that $f\circ g = g\circ f$, Result. Math. {\bf 31} (1997), 14-27. 
\bibitem{FolSzi} S. Foldes and J. Szigeti, Which self-maps appear as lattice anti-endomorphisms?, Algebra Univers. {\bf 75} (2016), 439-449.
\bibitem{Nat} M.B. Nathanson, \textit{Elementary Methods in Number Theory}, Springer, 2000. 
\bibitem{Newman_IM} I. Newman, \textit{Integral matrices}, Academic Press, 1972. 
\bibitem{Szigeti} J. Szigeti, Which self-maps appear as lattice endomorphisms?, Discrete Math.,  \textbf{321} (2014), 53-56.
\end{thebibliography}
\end{document}